\documentclass[12pt]{article}

\usepackage{graphicx}
\usepackage{amsmath,amsthm}
\usepackage{amssymb}

\newtheorem{thm}{Theorem}
\newtheorem{prop}{Proposition}
\newtheorem{alg}{Algorithm}
\newtheorem{cor}{Corollary}

\newcommand{\dk}[2]{\genfrac{\lbrack}{\rbrack}{0pt}{}{#1}{#2}}
\newcommand{\dl}[2]{\genfrac{\langle}{\rangle}{0pt}{}{#1}{#2}}
\newcommand{\PG}{\mbox{$\mathit{PG}(r-1,q)$}}
\newcommand{\AG}{\mbox{$\mathit{AG}(r-1,q)$}}
\newcommand{\set}[1]{\lbrace#1\rbrace}

\begin{document}

\title{Hyperplanes and hamiltonian circuits in Perfect Matroid Designs with fixed basis} 
\author{Wojciech Kordecki\\
University of Business in Wroc{\l}aw 
}

\maketitle

\begin{abstract}
We study the number of hamiltonian circuits, containing a fixed basis, and the number of hyperplanes, which do not contain a fixed basis in perfect matroid designs. 
Projective and affine finite geometries are considered as examples of such matroids. We give algorithms to find the hyperplanes and the hamiltonian circuits in such cases.
\end{abstract}

\section{Introduction}

In this note we use a standard matroid terminology -- see \cite{Oxley:ma2nd}
and \cite{Welsh:ma}.
Let $M=\left(E,\mathcal{F}\right)$ denote a matroid on a finite nonempty set $E$, where
$\mathcal{F}$ denotes a family of flats. By $\rho\left(A\right)$ we denote a rank of a
set $A$.
A family of circuits is denoted by $\mathcal{C}$, a family of bases is denoted
by $\mathcal{B}$, and  $\sigma\left(A\right)$ denote a span of $A$.
A flat of rank $\rho\left(E\right)-1$ is called a hyperplane.

We also need some definitions from the theory of projective
and affine finite geometries.
The monograph of Hirschfeld \cite{Hir:ProjGe-2nd} (see also \cite{Batten:finite_geo}) gives a detailed and
self-contained exposition of finite projective geometries.

Let $\mathit{GF}\left(q\right)$ be a Galois field, where $q$ is a prime power and let
$V\left(r,q\right)$ be an $r$-dimensional vector space on $\mathit{GF}\left(q\right)$. By $\mathcal{L}_P$ we
denote the lattice of subspaces of $V$.
Atoms of $\mathcal{L}_P$ constitute the points of  projective
geometry $\PG$ of dimension $r-1$.
Let $e$ be an element of the lattice $\mathcal{L}_P$. If $e$
is a subspace of dimension $k-1$ and $k>2$ we
define subspace of rank $k$ of $\PG$ as a set $A$ of all
points $p<e$, then $\rho\left(A\right)=k$.
Recall that the subspaces of $\PG$ form a modular
geometric lattice, hence if $A,B\in \mathcal{J}$, then $\rho(A\cup
B)=\rho\left(A\right)+\rho\left(B\right)-\rho\left(A\cap B\right)$, (see Welsh \cite{Welsh:ma}, p.~195).
The affine geometry $\AG$ is obtained from $\PG$ by deleting from the
latter all points of a hyperplane.
A line is a subspace of rank 2 and for
$x\neq y$, denote $l=L(x,y)=\sigma(\{x,y\})$. Note
that every 2-element set is independent.

Let
\[
[x]=\frac{q^x-1}{q-1}.
\]
Gaussian coefficients are defined as
\[
\dk xk =\prod_{j=1}^k \frac{q^{x-j+1}-1}
{q^j-1} \text{  for }0\le k\le x\,.
\]

It is well known that $\PG$ and $\AG$ are matroids $\left(E,\mathcal{F}\right)$ where
$E$ is
the set of points of the geometry and $\mathcal{F}$ is the family of its
subspaces.
It is also well known that the ground set $E$ of $\PG$ has $[r]$ elements, $[r]$
hyperplanes and $\dk rk$ subspaces of rank $k$. Similarly, $\AG$ has 
$q^{r-1}$ elements, $q^{r-k}\dk{r-1}{k-1}$ subspaces of rank $k$ and
$q[r-1]$ hyperplanes.

Let $S_1$ and $S_2$ are subspaces of $\PG$. By $S_1\vee
S_2=\sigma(S_1,S_2)$ we denote the smallest subspace containing both $S_1$
and $S_2$.
Let $B=\set{x_1,\dots,x_r}$ be any basis of $\PG$ and let $H_{\alpha}$ be a
hyperplane.

Kelly and Oxley give (\cite{KeOx:Asymp}, Lemma D), the following result.
\begin{prop}\label{prop:circall}
If $B$ is a basis of $\PG$, then there are precisely $\left(q-1\right)^{r-1}$
points $p$ of $\PG$ such that $B\cup\{p\}$ is a circuits.
\end{prop}
The following version of
Proposition \ref{prop:circall} for $\AG$ is given by Oxley in \cite{Oxley:ma2nd}, (Lemma~6.2.4, p.~172).
\begin{prop}\label{prop:circall-ag}
If $B$ is a basis of $\AG$ then there are precisely
$q^{-1}\left((q-1)^r-(-1)^r\right)$
points $p$ of $\AG$ such that $B\cup\{p\}$ is a circuits.
\end{prop}

\section{Results}
In this section we generalise Propositions \ref{prop:circall} and
\ref{prop:circall-ag} to the class of symmetric Perfect Matroid Designs, which includes
projective and affine geometries as particular cases.

A perfect matroid design or PMD is a matroid $M$ in which $k$-flat has a
common cardinal $\alpha_k$, $1\leq k\leq\rho(M)$. Such objects were first
studied by Young, Murthy and Edmonds
in \cite{YME:Equicardinals}, (see Welsh \cite{Welsh:ma}, p.~209 and Cameron and Deza \cite{ColbournDinitz:Handbook}, Chapter VII.10).

Let $M$ be a matroid of rank $r$, such that following
conditions hold:
\begin{itemize}
\item[(i)] all flats of the rank $j$ in $M$, $j<r$ are isomorphic,
\item[(ii)] the number of flats of rank $s$ in $M$ which contain a fixed flat
$M'$ of rank $u$ is equal $f\left(r,s,u\right)$ and does not depend on selection of
$M'$.
\end{itemize}
Call a matroid $M$ if it fulfills the above conditions, Symmetric Perfect Matroid Design. We denote it briefly by SPMD. Note, that both $\PG$ and
$\AG$ are SPMD.

Let $\dl{r}{k}$ denote the number of flats of rank $k$
(Whitney numbers of the second kind)
in a SPMD $M$ of
rank $r$. Let $\langle k\rangle=\dl{k}{1}$, denote the number of
elements of such $M$. Consider the number $f\left(r,s,u\right)$ of flats of rank $s$
which contains a fixed flat of rank $u$.
From the obvious equality
\[
\dl{r}{s}\dl{s}{u}=\dl{r}{u}f\left(r,s,u\right)
\]
we obtain the following formula:
\begin{equation}\label{eq:contain}
f\left(r,s,u\right)=\frac{\dl{r}{s}\dl{s}{u}}{\dl{r}{u}}\,.
\end{equation}
If $M$ is $\PG$ or $\AG$ and $u>0$, then from \eqref{eq:contain} we
get
\[f(r,s,u)=\dk{r-u}{s-u}\]
and for $\AG$ we have
\[
f\left(r,s,0\right)=q^{r-s}\dk{r-1}{s-1}\,.
\]
\begin{thm}\label{thm:hypergen}
If $B$ is a basis of SPMD $M$ of rank $r$, then there are
precisely
\begin{equation}\label{eq:hypergen}
n=\sum_{k=0}^{r-1}(-1)^k\frac{\binom{r}{k}\dl{r}{r-1}\dl{r-1}{k}}{\dl{r}{k}}
\end{equation}
hyperplanes $H$ which contain no elements of $B$.
\end{thm}
\begin{proof}
There are $\dl{r}{r-1}$ hyperplanes altogether
$rf(r,r-1,1)$ among them
contain at least one element from $B$, $\binom{r}{2}f(r,r-1,2)$
contain at least one pair of elements from $B$, and so on.
Hence from the inclusion-exclusion principle we get
\[
n=\sum_{k=0}^{r-1} \left(-1\right)^k\binom{r}{k}f\left(r,r-1,k\right)\,.
\]
Hence, from \eqref{eq:contain}, we obtain the assertion.
\end{proof}
\begin{thm}\label{thm:circgen}
If $B$ is a basis of SPMD $M$ of rank $r$ then there are
precisely
\begin{equation}\label{eq:circgen}
n=\sum_{k=0}^r(-1)^{r-k}\binom{r}{k}\left(\langle k\rangle-k\right)
\end{equation}
elements $e$ such that $B\cup e$ is a circuit.
\end{thm}
\begin{proof}
There are $\binom{r}{k}\left(\langle k\rangle-k\right)$ circuits, containing at
most $k$ elements from the basis $B$. Since $\langle1\rangle=1$ and
$\langle0\rangle=0$, we obtain the assertion.
\end{proof}
From Theorem \ref{thm:hypergen} we immediately obtain the following
result.
\begin{cor}\label{cor:hyperall}
If $B$ is a basis of $\PG$, then there are precisely $(q-1)^{r-1}$
hyperplanes $H$ which do not contain any element of $B$.
\end{cor}
\begin{proof}
For $\PG$ or $\AG$ we have
\begin{eqnarray}\label{eq:hyperall:pg}
n&=&\sum_{k=0}^{r-1} (-1)^k\binom{r}{k}\frac{q^{r-k}-1}{q-1} \nonumber \\
&=&\frac{1}{q-1}\left(\sum_{k=0}^{r-1} \left(-1\right)^k\binom{r}{k} q^{r-k}
-\sum_{k=0}^{r-1} \left(-1\right)^k\binom{r}{k}\right) \\
&=&(q-1)^{r-1}\,.  \nonumber
\end{eqnarray}
Hence, from \eqref{eq:hyperall:pg}, the assertion follows.
\end{proof}

Note, that  Proposition \ref{prop:circall} and Corollary \ref{cor:hyperall} are
strictly related. 
``Principle of duality'' (see Hirschfeld \cite{Hir:ProjGe-2nd}) establishes the one-to-one correspondence between points and hyperplanes
-- $p\leftrightarrow H$ and
$p_1\vee p_2\leftrightarrow H_1\cap H_2$.

Since $\AG$ is a SPMD matroid, then from
Theorem \ref{thm:hypergen} we have the following result.
\begin{cor}\label{cor:hyperall-ag}
If $B$ is a basis of $\AG$, then there are precisely
$n=\left(q-1\right)^{r-1}-1$ hyperplanes $H$ which contain no elements of $B$.
\end{cor}
\begin{proof}
From \eqref{eq:hyperall:pg} it follows that $\PG$ contains at least one
hyperplane which does not contain any element of $B$, (and it contains exactly
one if $q=2$). $\AG$ has the one such hyperplane less than $\PG$ has,
which establishes the formula.
\end{proof}
Arguing as above,
from Theorem \ref{thm:circgen} one can obtain Propositions
\ref{prop:circall} and \ref{prop:circall-ag} as well.
The proofs are similar to the proof of Proposition \ref{cor:hyperall}. Note,
that the proof of Proposition \ref{prop:circall-ag} given in \cite{Oxley:ma2nd}
is algebraic in nature, while proof of our more general
Theorem \ref{thm:circgen} has a simple combinatorial form.

\section{Algorithms for $\PG$ and $\AG$}
Until now all known geometric PMD are free matroids, $\PG$, $\AG$, Steiner system $S\left(t,k,v\right)$ and triffids of typ $\left(1,3,9,3^n\right)$ 
(see  Deza \cite{Deza:care}).
In this section we consider the cases of $\PG$ and $\AG$ (which seem the most interesting) in detail.
The following two algorithms give constructive methods to obtain hyperplanes described in Corollary \ref{cor:hyperall} and points described in Proposition \ref{prop:circall} in the case of $\PG$.
\begin{alg}\label{alg:hyperall:pg}~~  \\
{\bf Input:\/} Basis $B=\set{x_1,\dots,x_r}$ of $\PG$.\\
{\bf Output:\/} A family of all hyperplanes $\set{H_{\alpha}}$
in $\PG$ such that $\forall_{\alpha}:\,H_{\alpha}\cap B=\emptyset$.
\begin{enumerate}
\item Let $\left(x'_1,\dots,x'_r\right)$ be an arbitrary permutation of
$B$.
\item Let $l_i=L\left(x'_i,x'_{i+1}\right)$, $1\leq i\leq r-1$.
\item Choose $a_i=\left(y_{i,1},\dots,y_{i,q-1}\right)$, an arbitrary
sequence of points lying on $l_i$ different from $x'_i$ and $x'_{i+1}$,
$1\leq i\leq r-1$.
\item Choose ${\alpha}=(m_1,\dots,m_{r-1})$, an arbitrary sequence
such that $1\leq m_i\leq q-1$.
There are $\left(q-1\right)^{r-1}$ of such sequences.
\item Let $I_{\alpha}=\set{y_{1,m_1},\dots,y_{r-1,m_{r-1}}}$. The
set $I_{\alpha}$ is independent.
\item Return
$H_{\alpha}=\sigma\left(I_{\alpha}\right)=y_{1,m_1}\vee\dots\vee y_{r-1,m_{r-1}}$.
\end{enumerate}
\end{alg}
\begin{proof}
Since $y_{i,m_i}\in L\left(x'_i,x'_{i+1}\right)$, we have
\[
\set{y_{1,m_1},\dots,y_{k,m_k}}\subset\set{x'_1,\dots,x'_k}
\]
for $1\leq k\leq r-1$. Hence
$y_{k,m_k}\notin y_{1,m_1}\vee\dots \vee y_{k,m_k}$ for every
$k\leq r-1$.

Next we have to prove that for every $i$, $x'_i\notin H_{\alpha}$. Suppose,
that it is not the case, i.e.
there exists $x'_i\in H_{\alpha}$. Since from
step 3 and step 5 there exists $y\in L\left(x'_i,x'_{i+1}\right)$ (or $y\in
L(x'_i,x'_{i-1})$) and therefore $x'_{i+1}\in H_{\alpha}$, (or $x'_{i-1}\in
H_{\alpha}$). Hence all $x'_i\in H_{\alpha}$,
while we must have $\rho\left(H_{\alpha}\right)=r-1$.
\end{proof}

\begin{alg}\label{alg:circall:pg}~~ \\
{\bf Input:\/} Basis $B=\set{x_1,\dots,x_r}$ of $\PG$.\\
{\bf Output:\/} A family of all points $\{p_{\alpha}\}$ such that
$\forall_{\alpha}:\,B\vee \{p_{\alpha}\}$ is a circuit in $\PG$.
\begin{enumerate}
\item Let $\left(x'_1,\dots,x'_r\right)$ be an arbitrary permutation of
$B$.
\item Let $H'_i=\sigma\left(B\setminus x_i\right)$.
\item Choose $b_i=\left(H''_{i,1},\dots,H''_{i,q-1}\right)$, an arbitrary
sequence of $q-1$ hyperplanes, such that $H'_i\cap H'_{i+1}\subset
H''_{i,j}$ different from $H'_i$ and $H'_{i+1}$,
$1\leq i\leq r-1$.
\item Choose ${\alpha}=\left(m_1,\dots,m_{r-1}\right)$, an arbitrary sequence
such that $1\leq m_i\leq q-1$. There are $(q-1)^{r-1}$ of such sequences.
\item Let $p_{\alpha}=H''_{1,m_1}\cap\dots\cap H''_{r-1,m_{r-1}}$.
\item Return $p_{\alpha}$.
\end{enumerate}
\end{alg}
\begin{proof}
By the ``principle of duality'' and on the same way as in the
proof of Algorithm \ref{alg:hyperall:pg} we obtain the assertion.
\end{proof}

From both algorithms we can obtain the following important conclusion.
The families $\{H_{\alpha}\}$ and $\{p_{\alpha}\}$
produced by Algorithm \ref{alg:hyperall:pg} or
Algorithm \ref{alg:circall:pg} does not depend on the order of points selection
at the each step of algorithms.

Similar algorithms for $\AG$ can be constructed by a
straightforward modification of Algorithms \ref{alg:hyperall:pg}
and \ref{alg:circall:pg}. Since $\AG$ is obtained
from $\PG$ by deletion of a hyperplane, then we add a hyperplane to $\AG$
as a step~0. Next, after performing modified the last step~6, the hyperplane $H_0$ will be
removed.

\begin{alg}\label{alg:hyperall:ag} ~~ \\
{\bf Input:\/} Basis $B=\set{x_1,\dots,x_r}$ of $\AG$.\\
{\bf Output:\/} A family of all hyperplanes $\set{H_{\alpha}}$ in $\AG$ such
that $\forall_{\alpha}:\,H_{\alpha}\cap B=\emptyset$.
\begin{enumerate}
\setcounter{enumi}{-1}
\item Add $[r-1]$ points, which form a hyperplane $H_0$ in
$\PG$. 
\end{enumerate}
In the next steps, the $\PG$ arisen from $\AG$, is considered.
Steps 1 -- 5 are the same as in Algorithm~\ref{alg:hyperall:pg}.
\begin{enumerate}
\setcounter{enumi}{5}
\item If $H_{\alpha}\neq H_0$, then return
$H_{\alpha}=\sigma\left(I_{\alpha}\right)\setminus H_0$
\end{enumerate}
\end{alg}

\begin{alg}\label{alg:circall:ag} ~~ \\
{\bf Input:\/} Basis $B=\set{x_1,\dots,x_r}$ of $\AG$.\\
{\bf Output:\/} A family of all points $\{p_{\alpha}\}$ such that
$\forall_{\alpha}:\,B\vee \{p_{\alpha}\}$ is a circuit in $\AG$.
\begin{enumerate}
\setcounter{enumi}{-1}
\item Add $[r-1]$ points, which form a hyperplane $H_0$ in
$\PG$. 
\end{enumerate}
In the next steps, the $\PG$ arisen from $\AG$ is considered.
Steps 1 -- 5 are the same as in Algorithm~\ref{alg:circall:pg}.
\begin{enumerate}
\setcounter{enumi}{5}
\item If $p_{\alpha}\notin H_0$, then return $p_{\alpha}$.
\end{enumerate}
\end{alg}

\bibliography{hyper}
\bibliographystyle{plain}
Wojciech Kordecki \\
University of Business in Wroc{\l}aw \\
ul. Ostrowskiego 22 \\
53-238 Wroc{\l}aw \\
e-mail: wojciech.kordecki@handlowa.eu

\end{document}